\documentclass[12pt]{amsart}

\usepackage[margin=1.15in]{geometry}
\usepackage{amscd,amssymb, amsmath,amsfonts, wasysym, mathrsfs, mathtools,enumerate,hhline,color}

\usepackage{url}

\definecolor{hot}{RGB}{65,105,225}

\usepackage[pagebackref=true,colorlinks=true, linkcolor=hot ,  citecolor=hot, urlcolor=hot]{hyperref}

\theoremstyle{plain}
\newtheorem{theorem}{Theorem}[section]

\newtheorem{cor}[theorem]{Corollary}
\newtheorem{conj}[theorem]{Conjecture}
\newtheorem{lemma}[theorem]{Lemma}
\newtheorem{thrm}[theorem]{Theorem}

\theoremstyle{definition}

\newtheorem{defn}[theorem]{Definition}
\newtheorem{que}[theorem]{Question}
\newtheorem{rmk}[theorem]{Remark}

\newtheorem{exa}[theorem]{Example}
\newtheorem*{ex*}{Example}

\def\pa{\partial}

\newcommand\sO{{\mathcal O}}

\newcommand\sH{{\mathcal H}}

\newcommand\sM{{\mathcal M}}

\newcommand\sV{{\mathcal V}}

\newcommand\sR{\mathcal{R}}

\def\sQ{\mathcal{Q}}
\def\bR{\mathbf{R}}

\newcommand\qq{{\mathbb{Q}}}
\newcommand\zz{{\mathbb{Z}}}
\newcommand\rr{{\mathbb{R}}}
\newcommand\cc{{\mathbb{C}}}

\newcommand\nn{{\mathbb{N}}}

\def\bH{\mathbb{H}}

\newcommand{\mb}{\mathcal{M}_{\textrm{B}}}

\newcommand{\mdol}{\mathcal{M_{\textrm{Dol}}}}
\newcommand{\Art}{\mathsf{ART}}
\newcommand{\Set}{\mathsf{SET}}
\def\Exp{{\rm{Exp}}}
\def\al{\alpha}
\DeclareMathOperator{\codim}{codim}              
\DeclareMathOperator{\id}{id}                    

\DeclareMathOperator{\ad}{ad}

\DeclareMathOperator{\homo}{Hom}
\def\Hom{\homo}
\DeclareMathOperator{\enmo}{\mathcal{E}\hspace{-2pt}\it{nd}}

\DeclareMathOperator{\rank}{Rank}

\def\ra{\rightarrow}
\def\cU{\mathcal{U}}
\def\bone{\mathbf{1}}
\def\bC{\mathbb{C}}

\def\cM{\mathcal{M}}
\def\Def{{\rm {Def}}}
\def\cR{\mathcal{R}}

\def\lam{\lambda}
\def\lra{\longrightarrow}

\def\cO{\mathcal{O}}
\def\Pic{{\rm{Pic}}}
\def\cK{\mathcal{K}}
\def\cF{\mathcal{F}}
\def\Tor{{\rm{Tor}}}
\def\cV{\mathcal{V}}
\def\cA{\mathcal A}
\def\bA{\mathbb{A}}

\newcommand{\ubul}{{\,\begin{picture}(-1,1)(-1,-3)\circle*{2}\end{picture}\ }}

\newcommand{\Ap}{{{A}}^{p,\ubul}_{\rm{Dol}}}
\newcommand{\Ad}{{A}^{0,\ubul}_{\rm{Dol}}}
\newcommand{\Ar}{{A}^{\ubul}_{\rm{DR}}}

\def\bZ{\mathbb{Z}}

\def\ed{{\rm{exdeg}}}
\def\gr{{\rm{gr}}}
\def\cD{\mathcal{D}}

\title[Recent results on cohomology jump loci]{Recent results on cohomology jump loci}
\begin{document}

\author{Nero Budur}
\address{KU Leuven, Department of Mathematics,
Celestijnenlaan 200B, B-3001 Leuven, Belgium} 
\email{Nero.Budur@wis.kuleuven.be}

\author{Botong Wang}
\address{KU Leuven, Department of Mathematics,
Celestijnenlaan 200B, B-3001 Leuven, Belgium} 
\email{wang\_botong@hotmail.com}

\date{}

\begin{abstract} 
This is a survey of  recent results related to cohomology jump loci. It emphasizes connections with deformations with cohomology constraints, global structural results for rank one local systems and line bundles, some connections with restrictions on fundamental groups and homotopy types, and connections with classical singularity theory and Bernstein-Sato ideals.
\end{abstract}
\maketitle

\begin{center}
\emph{Dedicated to Steve Zucker's 65th birthday}
\end{center}

\setcounter{tocdepth}{1}
\tableofcontents

\section{Introduction}

Let $\cM$ be a moduli space parametrizing objects $P$ that come with $\bC$ vector spaces 
$$
H^i (P)\quad (i\in\bZ).
$$
For example, some cohomology theory. The {\bf cohomology jump loci} of $\cM$ are defined set-theoretically as
$$
\cV^i_k =\{ P\in \cM\mid \dim H^i(P)\ge k\}.
$$
For every $i$, there is a chain of inclusions
$$
\cM = \cV^i_0\supset \cV^i_1\supset\cV^i_2\supset\ldots
$$

In many situations, the cohomology jump loci form closed subsets of the moduli space. When the moduli is an algebraic variety or a scheme, the cohomology jump loci are usually closed subschemes. Cohomology jump loci are at least as complex as the moduli, and in general everything can be a moduli space. Hence the study of the cohomology jump loci can be described as finding a balance between general setups one can focus on and interesting properties one can draw.

In this article one finds a survey of the recent results related to cohomology jump loci. It emphasizes connections with deformations with cohomology constraints, global structural results for rank one local systems and line bundles, some connections with restrictions on fundamental groups and homotopy types, and connections with classical singularity theory and Bernstein-Sato ideals.

There are many results left out. For example, we do not survey the all the implications on quasi-projective groups one can draw from the global structure of cohomology jump loci of rank one local systems. 
We also do not survey the large body of work dedicated to complements of hyperplane arrangements. One can find these two topics in the surveys by Suciu \cite{Su,Su2}.

{\bf Acknowledgement.} The first author thanks Steve Zucker and the Department of Mathematics of the Johns Hopkins University for the hospitality during the writing of this survey.  The first author was partly sponsored by NSA, FWO, a KU Leuven OT grant, and a Flemish Methusalem grant.

\section{Cohomology jump ideals}\label{secIntro}

\subsection{Old algebraic example: Fitting ideals} Let $R$ be a commutative ring and $M$ a finitely generated $R$-module. Consider a free finite rank presentation of $M$
$$
R^n\mathop{\lra}^{d} R^m\lra M\lra 0.
$$
Define the {\bf Fitting ideals} of $M$ 
$$
J_k(M) \subset R
$$
as the ideals
$$
J_k(M) = I_{m-k+1}(d)
$$ 
generated by the minors of size $m-k+1$ of the matrix $d$. We set $J_k(M)=0$ for $k\le 0$, and $J_k(M)=R$ for $k> m$. One has the classical result of Alexander and Fitting, see \cite[\S 20.2]{E}:

\begin{thrm} 
Let $R$ be a commutative ring and $M$ a finitely generated $R$-module. The Fitting ideal $J_k(M)$ depends only on $M$ and $k$, and not on the presentation of $M$.
\end{thrm}

Let $\cM=Spec(R)$. The ideals $J_k(M)$ define closed subschemes 
$$\cV_k(M)\subset \cM.$$
There is a chain of closed embeddings
$$
\cM=\cV_0(M)\supset \cV_1(M)\supset\ldots\supset \cV_{m+1}(M)
$$
Thinking of $\cM$ as the moduli space of its closed points, each closed point comes with an assigned vector space
$$
(m=\text{maximal ideal of }R) \mapsto M\otimes _R R/m.
$$
Moreover,  $\cV_k(M)$ is set-theoretically
$$\cV_k(M)=\{ m\subset R \text{ maximal ideal }\mid \dim_{R/m}(M\otimes _R R/m)\ge k\},$$
hence it is a cohomology jump locus.

\subsection{Cohomology jump ideals are generalized Fitting ideals.}\label{subCJI} In \cite{BW2}, the construction and the key properties of the Fitting ideals were generalized simultaneously in two directions: to define invariants of a module from a whole free resolution, not just from a free presentation, and to define these invariants for complexes of modules.

Let $R$ be a noetherian commutative ring, and let $M^\ubul$ be a complex of $R$-modules, bounded above. Suppose $H^i(M^\ubul)$ is a finitely generated $R$-module. There always exists a bounded above complex $F^\ubul$ of finitely generated free $R$-modules and a morphism of complexes $$g: F^\ubul\mathop{\lra}^{\sim} M^\ubul$$ which is a quasi-isomorphism. We define the {\bf cohomology jump ideals} of $M^\ubul$ to be 
$$
J^i_k(M^\ubul)=I_{rank(F^i)-k+1}(d^{i-1}\oplus d^{i})
$$
where $I$ denotes the determinantal ideal as above, and  $$d^{i-1}: F^{i-1}\to F^i\quad\text{ and }\quad d^i: F^i\to F^{i+1}$$ are the differentials of the complex $F^\ubul$. The same conventions for ideals generated by minors as in the previous subsection apply when $rank(F^i)-k+1$ falls outside the range of the size of minors. 

The {\bf cohomology jump loci} of a bounded above complex $M^\ubul$ of $R$-modules with finitely generated cohomology are the closed subschemes
$$
\cV^i_k(M^\ubul)\subset Spec (R)
$$
defined by the cohomology jump ideals $J^i_k(M^\ubul)$. In some sense, this is the mother of all examples of cohomology jump loci.

\begin{rmk} If $M$ is a finitely generated $R$-module, then the Fitting ideals of $M$ coincide with the $0$-cohomology jump ideals of $M$ seen as a complex 
$$M[0]=[\ldots\ra 0\ra M\ra 0\ra\ldots]$$
with trivial terms outside degree $0$:
$$
J_k(M)=J^0_k(M[0])\quad\text{ and }\quad \cV_k(M)=\cV^0_k(M[0]).
$$ 
\end{rmk}

\begin{thrm}{\rm (\cite{BW2})} Let $R$ be a noetherian commutative ring. Let $M^\ubul$ be a complex of $R$-modules, bounded above, with finitely generated cohomology.  Then $J^i_k(M^\ubul)$ depends only on $M$, $i$, and $k$, and not on the choice of resolution $F^\ubul$ of $M^\ubul$. 
\end{thrm}

An important first property is:

\begin{cor}\label{idealindep}
If $M^\ubul$ is quasi-isomorphic to $N^\ubul$, then $J^i_k(M^\ubul)=J^i_k(N^\ubul)$. 
\end{cor}

Since taking determinantal ideals commutes with base change, one has the following base change property of the cohomology jump ideals:

\begin{cor}\label{tensor}
Let $R$ and $M^\ubul$ be defined as above, and let $S$ be a noetherian $R$-algebra. Moreover, suppose $M^\ubul$ is a complex of flat $R$-modules, then $$J^i_k(M^\ubul)\cdot S=J^i_k(M^\ubul\otimes_R S),$$ where we regard $M^\ubul\otimes_R S$ as a complex of $S$ modules. 
\end{cor}

When $R$ is a field, by definition $J^i_k(M^\ubul)=0$ if $\dim H^i(M^\ubul)\geq k$ and $J^i_k(M^\ubul)=R$ if $\dim H^i(M^\ubul)<k$. Thus, we have the following:

\begin{cor}\label{corFl}
Suppose $M^\ubul$ is a complex of flat $R$-modules. Then for any maximal ideal $m$ of $R$, $J^i_k(M^\ubul)\subset m$ if and only if $$\dim_{R/m} H^i(M^\ubul\otimes_R R/m)\geq k.$$ 
\end{cor}

This says that the cohomology jump ideals deserve their name.

\begin{rmk} (Tor jumping loci)
If $M$ is a finitely generated $R$-module, the ideals $J^i_k(M)$ put a closed subscheme structure on the Tor jumping loci:
$$
\cV^i_k(M)_{red}=\{  P\in Spec (R) \text{ closed point }\mid \dim_{R/m_P} \Tor_{-i}^R(M,R/m_P)\ge k \}.
$$
\end{rmk}

\begin{rmk} (Syzygies) If $M$ is a finitely generated $R$-module and $F^\ubul$ is a free resolution, the syzygies of $M$ are defined as the kernels of the differentials of $F^\ubul$. Hence the cohomology jump loci of a syzygy of $M$ coincide with the far-enough-to-the-left cohomology jump loci of $M$. 
\end{rmk}

\subsection{Dimension bounds.}\label{subDBH} Recall the following classical homological algebra statement from  \cite[Theorem 20.9]{E}:

\begin{thrm} 
Let $R$ be a ring. Let 
$$
0\lra F^{0} \mathop{\lra}^{d^{0}} F^{1}\lra \ldots \lra F^{n-1} \mathop{\lra}^{d^{n-1}} F^{n}
$$
be a complex of free finitely generated $R$-modules. Then this complex is exact if and only if
$$
rank (F^i)=rank (d^{i-1}\oplus d^i)
$$
and
$$
depth \;I_{rank(d^i)}(d^i)\ge n-i.
$$
\end{thrm}

By applying this statement to the exact complex $(F^\ubul[-1]\oplus F^\ubul, d[-1]\oplus d)$, one has the following universal bound for cohomology support loci:

\begin{thrm}\label{thmDim}
With the same assumptions as in the previous theorem, if $F^\ubul$ is exact, then
$$
depth\; J^i_1(F^\ubul)\ge n-i.
$$
\end{thrm}

\section{Deformation theory}

\subsection{Local study.}\label{subDef} By a {\bf deformation problem} we mean describing the formal scheme or the analytic germ $\cM_{(P)}$ at some object $P$ of some moduli space $\cM$. A deformation problem is equivalent to describing the corresponding functor, which we denote also by $\cM_{(P)}$:
$$
\cM_{(P)} : \Art \lra \Set,
$$
where $\Art$ is the category of Artinian local algebras of finite type over $k$, and $\Set$ is the category of sets. Assuming $\cM$ is a $k$-scheme, on objects the functor is defined as
$$
A\mapsto \homo_{k-{\text{sch}}}(Spec(A),\cM_{(P)})=\homo _{k-\text{alg}}(\cO_{\cM_{(P)}},A).
$$
This functor is usually well-defined even if the moduli space is not. That is, by passing to the categorical language, one can still ask about the deformations of some object $P$ inside a category. For example, one can ask about the deformations of a coherent sheaf on an algebraic variety, regardless of any stability conditions necessary for the construction of moduli spaces.

Suppose now the objects parametrized by $\cM$ come with a cohomology theory, so that one has the stratification of $\cM$ by the cohomology jump loci $\cV^i_k$. The study of the formal or analytic germs of the cohomology jump loci at an object $P$ amounts to studying the associated functor on Artinian local algebras as a subfunctor  $$
\cV^i_{k,(P)}\subset \cM_{(P)}.$$ Even if the moduli space does not exist, the corresponding deformation problem with cohomology constraints is usually well-defined.

\subsection{Differential graded Lie algebras.}\label{subsDGLA} In a letter written by P. Deligne to J. Millson in 1986 in connection with the paper \cite{gm} by Goldman-Millson, a principle of sweeping generality was stated. The letter is available on J. Millson's webpage \cite{Del}.

\medskip
\noindent
{\bf Deligne's Principle. }{\it Every deformation problem over a field of characteristic zero is controlled by a differential graded Lie algebra, with two quasi-isomorphic differential graded Lie algebras giving the same deformation theory.}
\medskip

 We explain next the meaning of Deligne's Principle.  
 
 From now on, $k=\bC$. We use the abbreviation DGLA for a {\bf differential graded Lie algebra}. 
\begin{defn}
A {\bf DGLA} consists of the following set of data,
\begin{enumerate}
\item a graded vector space $C=\bigoplus_{i\in \mathbb{N}}C^i$ over $\mathbb{C}$,
\item a Lie bracket which is bilinear, graded skew-commutative, and satisfies the graded Jacobi identity, i.e., for any $\alpha\in C^i, \beta\in C^j$ and $\gamma\in C^k$,
$$[\alpha, \beta]+(-1)^{ij}[\beta, \alpha]=0$$
and
$$(-1)^{ki}[\alpha, [\beta, \gamma]]+(-1)^{ij}[\beta, [\gamma, \alpha]]+(-1)^{jk}[\gamma, [\alpha, \beta]]=0$$
\item a family of linear maps, called the differential maps, $d^i: C^i\to C^{i+1}$, satisfying $d^{i+1}d^i=0$ and the Leibniz rule, i.e., for $\alpha\in C^i$ and $\beta\in C$
$$d[\alpha, \beta]=[d\alpha, \beta]+(-1)^i[\alpha, d\beta]$$
where $d=\sum d^i: C\to C$. 
\end{enumerate} We denote this DGLA by $(C, d)$, or $C$ when there is no risk of confusion. 
A homomorphism of DGLAs is a linear map which preserves the grading, Lie bracket, and the differential maps. 
\end{defn}

\begin{defn}
A homomorphism of DGLAs $g: C\to D$ is {\bf 1-equivalent} if it induces an isomorphism on cohomology up to degree $1$ and a monomorphism at degree $2$. Two DGLAs are of the same $1$-\textbf{homotopy type}, if they can be connected by a zig-zag of $1$-equivalences. 
\end{defn}

The nicest DGLAs form a particular class:

\begin{defn}
For a DGLA $(C,d)$, the cohomology forms a natural DGLA $(H^\ubul(C),0)$ with zero differential. We say the DGLA $C$ is $1$-\textbf{formal} if $C$ is of the same $1$-homotopy type as $(H^\ubul(C), 0)$. 
\end{defn}

Given a DGLA  $(C,d)$, one can attach naturally a functor, called the {\bf deformation functor},
$$
\Def(C,d):\Art\lra \Set
$$
defined on Artinian local rings by
$$
A\mapsto  \frac{\text{Maurer-Cartan elements of }C\otimes_\bC m_A}{\text{gauge}},$$
that is,
$$
A\mapsto \frac{\{ \omega\in C^1\otimes_{\cc}m_A\;|\; d\omega+\frac{1}{2}[\omega, \omega]=0\}}{C^0\otimes_\bC m_A}.
$$
Here $m_A$ is the maximal ideal of $A$, and $C\otimes_{\cc}m_A$ is naturally a DGLA by letting 
\begin{align*}
[\alpha\otimes a, \beta\otimes b]&=[\alpha, \beta]\otimes ab\\
d(\alpha\otimes a)&=d\alpha \otimes a.
\end{align*}
Since $(C\otimes_{\cc}m)^0=C^0\otimes_{\cc}m$ is a nilpotent Lie algebra, the Campbell-Hausdorff multiplication defines a nilpotent Lie group structure on the space $C^0\otimes m$. The gauge action of an element $\lambda\in C^0\otimes m$ on $C^1\otimes m$ is given by
$$\alpha \mapsto \exp(\ad \lambda)\alpha+\frac{1-\exp(\ad \lambda)}{\ad \lambda}(d\lambda)$$
in terms of power series.

Suppose now one has a deformation problem $\cM_{(P)}$ over $\bC$ as in \ref{subDef}. To say that a DGLA $C$ controls this deformation problem means that the two functors 
$$
\cM_{(P)} ) :\Art \to\Set\quad\quad\text{and}\quad\quad \Def(C):\Art \to\Set
$$ 
are isomorphic:
$$
\cM_{(P)} = \Def(C).$$

Constructing a DGLA responsible for a deformation problem is an ad-hoc process with no all-encompassing recipe. Typically the DGLAs one first constructs in deformation theory are infinite-dimensional, and there might be multiple choices of such constructions. The second statement of Deligne's Principle is then crucial:

\begin{theorem}[\cite{gm}]\label{gm0}
The deformation functor $\Def(C)$ only depends on the 1-homotopy type of $C$. More precisely, if a morphism of DGLA $f: C\to D$ is 1-equivalent, then the induced transformation on functors $f_*: \Def(C)\to \Def(D)$ is an isomorphism. 
\end{theorem}

This finishes the explanation of Deligne's Principle.

\begin{exa}\label{exFIM} (\cite{FIM}) For a complex manifold or complex algebraic variety $X$, the infinitesimal deformations of an $\cO_X$-coherent sheaf $E$ are controlled by the DGLA  of global sections $$\Gamma (X,\cA^\ubul (\enmo^\ubul (\tilde{E}^\ubul)))$$ of any acyclic resolution $\cA^\ubul$ of the sheaf of DGLAs $\enmo^\ubul(\tilde{E}^\ubul)$ of a locally free resolution  $\tilde{E}^\ubul$ of $E$. If $X$ is a complex manifold, then $\cA^\ubul$ can chosen to be the Dolbeault resolution. 
\end{exa}

\begin{rmk}
A theorem due to  Lurie \cite{Lu} and Pridham \cite{Pr} in the framework of derived algebraic geometry states that, with the appropriate axiomatization, Deligne's Principle is a theorem. However, in this axiomatization the notion of a deformation problem might be different than the classical notion we have defined, and a classical deformation problem might not come from any new type of deformation problem.
\end{rmk}

\subsection{Deformations  with constraints.}\label{subCC}  By a {\bf deformation problem with cohomology constraints} we mean the study of formal germs of cohomology jump loci
$$\sV^i_{k ,(P)}\subset \cM_{(P)}$$
for all $i$ and $k$ at once. This is equivalent to describing the corresponding functors on local Artin rings which we denote also by $\sV^i_{k ,(P)}$, and which in fact are usually well-defined even if the formal germs are not. 

In \cite{BW2}, a generalization of Deligne's Principle has been proposed and illustrated in several examples:

\medskip
\noindent
{\bf New Principle. }{\it Every deformation problem with cohomology constraints over a field of characteristic zero is controlled by a pair consisting of a DGLA together with a DGL module, with two quasi-isomorphic pairs giving the same deformation theory with cohomology constraints.}
\medskip

Let us explain what this means. The deformation problem without constraints $\cM_{(P)}$ comes, as claimed by Deligne's Principle, with a DGLA $C$. Typically this is some kind of endomorphism object. That is, $C$ does not come alone, but it comes with a natural module $M$ on which it acts. Given a pair $(C,M)$ consisting of a DGLA and  DGL module, \cite{BW2} defines canonicallly {\bf cohomology jump deformation functors} 
$$\Def^i_k(C,M):\Art\to\Set .$$ 
Then $$\sV^i_{k ,(P)}=\Def^i_k(C,M)\quad\text{ for all }i, k$$ as subfunctors of $\cM_{(P)}=\Def(C)$. We recall now the basic details.

\subsection{Deformation theory via DGLA pairs.}\label{subDefPairs} 

\begin{defn}\label{module}
Given a DGLA $(C, d_C)$, we define a \textbf{module} over $(C, d_C)$ to be the following set of data,
\begin{enumerate}
\item a graded vector space $M=\bigoplus_{i\in \mathbb{N}} M^i$ together with a bilinear multiplication map $C\times M\to M$, $(a, \xi)\mapsto a\xi$, such that for any $\alpha\in C^i$ and $\xi\in M^j$, $\alpha\xi \in M^{i+j}$. And furthermore, for any $\alpha\in C^i, \beta\in C^j$ and $\zeta\in M$, we require
$$[\alpha, \beta]\zeta=\alpha(\beta\zeta)-(-1)^{ij}\beta(\alpha\zeta).$$
\item a family of linear maps $d^i_M: M^i\to M^{i+1}$ (write $d_M=\sum_{i\in\zz} d^i_M: M\to M$), satisfying $d^{i+1}_M d^i_M=0$. And we require it to be compatible with the differential on $C$, i.e., for any $\alpha\in C^i$, 
$$d_M(\alpha\xi)=(d_C\alpha)\xi+(-1)^i\alpha(d_M\xi).$$
\end{enumerate}
\end{defn}

\begin{defn}
A \textbf{homomorphism} of $(C, d_C)$-modules $f: (M, d_M)\to (N, d_N)$ is a linear map $f: M\to N$ such that
\begin{enumerate}
\item $f$ preserves the grading, i.e., $f(M^i)\subset N^i$,
\item $f$ is compatible with multiplication by elements in $C$, i.e., $f(\alpha\xi)=\alpha f(\xi)$, for any $\alpha\in C$ and $\xi \in M$,
\item $f$ is compatible with the differentials, i.e., $f(d_M\alpha)=d_N f(\alpha)$.
\end{enumerate}
\end{defn}

Fixing a DGLA $(C, d_C)$, the category of $C$-modules is an abelian category. 

\begin{defn}
A \textbf{DGLA pair} is a DGLA $(C, d_C)$ together with a $(C, d_C)$-module $(M, d_M)$. Usually, we denote such a pair simply by $(C, M)$.  A homomorphism of DGLA pairs $g: (C, M)\to (D, N)$ consists of a map $g_1: C\to D$ of DGLA and a $C$-module homomorphism $g_2: M\to N$, considering $N$ as a $C$-module induced by $g_1$. For $q\in \nn\cup \{\infty\}$, we call $g$ a $q$-\textbf{equivalence} if $g_1$ is 1-equivalent and $g_2$ is $q$-equivalent. Here a morphism of complexes is $q$-equivalent if it induces an isomorphism on cohomology up to degree $q$ and a monomorphism at degree $q+1$ (for example, $\infty$-equivalent means quasi-isomorphic, which also means homotopy equivalence as complexes since these are complexes of vector spaces). Moreover, we define two DGLA pairs to be of the same $q$-\textbf{homotopy type}, if they can be connected by a zig-zag of $q$-equivalences. Two DGLA pairs have the same {\textbf{homotopy type}} if they have the same $\infty$-homotopy type.
\end{defn}

By keeping the 1-equivalence for the DGLAs in the definition of the $q$-equivalence for DGLA pairs, we are using in an economical manner the fact that the Maurer-Cartan set (to be introduced below) of a DGLA is independent under 1-equivalences. 

\begin{defn} Let $(C, M)$ be a DGLA pair. Then $(H^\ubul(C),d=0)$, the cohomology of $C$ with zero differentials, is a DGLA, and $(H^\ubul(M),d=0)$, the cohomology of $M$ with zero differentials, is an $H^\ubul(C)$-module. We call the DGLA pair $$(H^\ubul(C), H^\ubul(M))$$ the \textbf{cohomology DGLA pair} of $(C, M)$. 
\end{defn}

\begin{defn}
We say the DGLA pair $(C, M)$ is $q$-\textbf{formal} if $(C, M)$ is of the same $q$-homotopy type as $(H^\ubul(C), H^\ubul(M))$. A pair is {\bf formal} is it is $\infty$-formal.
\end{defn}

From now, for a DGLA pair $(C,M)$ we always assume that $M$ is bounded above as a complex and $H^j(M)$ is a finite dimensional $\cc$-vector space for every $j\in \zz$. This is needed in order to apply the results about cohomology jump ideals from \ref{subCJI}.

\begin{lemma} Let $(C,M)$ be a DGLA pair. Let $A\in\Art$. Given any Maurer-Cartan element $\omega$ of $C\otimes_\bC m_A$, that is
$$
\omega\in C^1\otimes_\bC m_A\text{ with } d\omega +\frac{1}{2}[\omega,\omega]=0,
$$
there is a complex of $A$-modules 
\begin{equation}\label{eqAo}
(M\otimes_\bC A, d_\omega)
\end{equation}
with  $$d_\omega:=d\otimes \id_A+\omega.$$ Moreover, $(M\otimes_\bC A,d_\omega)$ has finitely generated cohomology over $A$.
\end{lemma}

In particular, one can apply the construction of cohomology jump ideals from \ref{subCJI} to obtain ideals
$$
J^i_k(M\otimes_\bC A, d_\omega) \subset A.
$$
The {\bf cohomology jump deformation functors} of the DGLA pair $(C,M)$ 
$$
\Def^i_k(C,M): \Art \to \Set
$$
are defined by
$$
A\mapsto \frac{\{ \omega\in C^1\otimes_{\cc}m_A\;|\; d\omega+\frac{1}{2}[\omega, \omega]=0\text{ and } J^i_k(M\otimes_\bC A, d_\omega)=0\}}{C^0\otimes_\bC m_A}.
$$
This gives a well-defined subfunctor of $\Def(C)$, see \cite{BW2}. Theorem \ref{gm0} admits the following refinement:
 
\begin{theorem}\label{thmHot}
The cohomology jump functor $\Def^i_k(C, M)$ only depends on the number $k$ and on the $i$-homotopy type of $(C, M)$. \end{theorem}

\subsection{Finite-dimensional DGLA pairs.}\label{subEqs} In the presence of finite dimensional DGLA pairs,  the cohomology jump deformation functors are pro-represented by schemes with explicit equations, up to the gauge action. We explain what this means.

Let $(C,M)$ be a finite dimensional DGLA pair. Fixing a base $e_1,\ldots,e_b$ of $C^1$, denote by $x_1,\ldots, x_b$ the dual base. We identify the affine coordinate ring $\cO_{C^1}$ with the polynomial ring $\bC[x_1,\ldots,x_b]$. The universal element of $C^1$ is defined as the element
$$
w_{univ}=\sum_{i=1}^be_i\otimes x_i \quad\text{in}\quad C^1\otimes_\bC\cO_{C_1}.
$$
The universal Maurer-Cartan equations are obtained by setting in $C^2\otimes_\bC\cO_{C^1}$,
$$
(d_{C}\otimes id_{\cO_{C^1}})w_{univ}+\frac{1}{2}[w_{univ},w_{univ}]=0,
$$ 
where the Lie bracket is extended by multiplication to $C\otimes _\bC\cO_{C^1}$. The universal Maurer-Cartan equations define a closed affine subscheme
$$
\cF(C)\subset C^1\cong \bA_\bC^{b}
$$
containing the origin. This is also called the {\it versal deformation, or Kuranishi space} of $C$. 

$\Def(C)$ is pro-represented by $\cF(C)_{(0)}$ modulo the gauge action when this quotient exists in the category of schemes. This issue does not arise when, for example, the gauge action of $C^0$ on $C^1$ is trivial. In that case,
$$
\Def(C)\cong \cF(C)_{(0)}.
$$

If $C$ is 1-formal, then the universal Maurer-Cartan equations become the quadratic (that is, degree at most two) equations
$$
[w_{univ},w_{univ}]=0.
$$
In this case, $\cF(C)$ is a quadratic cone and we denote it by $\sQ(C)$.

On $\cF(C)$, one can define the $\cO_{\cF(C)}$-linear {\it universal Aomoto complex}
$$
(M^\ubul \otimes_\bC \cO_{\cF(C)},  d_M\otimes id + w_{univ}).
$$
Taking the cohomology jump ideals $$J^i_k(C,M)\subset\cO_{\cF(C)}$$ of this complex, one obtains conic affine closed subschemes
$$
\cF^i_k(C,M)\subset \cF(C).
$$
Now, the cohomology jump deformation functors $\Def^i_k(C,M)$ are pro-represented by $\cF^i_k(C,M)$ modulo the gauge action, when this quotient exists as a scheme. When the action of $C^0$ on $C^1$ is trivial, 
$$
\Def^i_k(C,M)\cong \cF^i_k(C,M)_{(0)}.
$$

When the DGLA pair is formal, the maps in the universal Aomoto complex are given just by multiplication with the image of $w_{univ}$ in $\cO_{\sQ(C)}$. Hence the entries of the matrices representing these maps are images of linear forms from $\cO_{C^1}=\bC[x_1,\ldots, x_b]$. The $J^i_k(C,M)$ are generated by minors of these matrices, as explained in \ref{secIntro}. We denote the loci $\cF^i_k(C,M)$ in this case by $$\cR^i_k(C,M),$$
and they are classically called {\it resonance varieties}.

\begin{rmk}
One can define the quadratic cones and the resonance varieties from cohomology in the non-formal cases as well, but the web of connections with original problem is in general largely lost and their geometric meaning is harder to grasp.
\end{rmk}

\subsection{Dimension bounds.}\label{subDB} We address the question of providing bounds for the dimensions of  cohomology jump loci $$\sV^i_k\subset \cM$$
locally around the point $P$.  Assume for simplicity that the controlling DGLA pair $(C,M)$ for these loci is a formal DGLA pair with trivial action of $H^0C$ on $H^1C$. Recall from the last subsection that it is enough to bound the dimension the resonance varieties of the cohomology pair,
$$
\sR^i_k=\sR^i_k(HC,HM)\subset \sQ(HC)=\sQ.
$$
Define the {\bf exactness degree} $$\ed(C,M)$$ to be the smallest degree at which the universal Aomoto complex for the cohomology pair
$$
(H^\ubul M \otimes_\bC \cO_{\sQ},   w_{univ}\cdot)
$$
 is not exact. By Theorem \ref{thmDim}, the depth of the ideal defining $\sR^i_1$ in $\sQ$ is $\ge \ed(C,M)-i$. Suppose in addition that $\sQ$ is nice enough so that the depth of ideals in $\cO_{\sQ}$ is the same as codimension, for example, that $\cO_{\sQ}$ is Cohen-Macaulay. Then locally at $P$ one has
\begin{equation}\label{eqN}
\codim \sV^i_{1} \ge \ed(C,M)-i.
\end{equation}

\section{Deformations of holomorphic vector bundles}

Let $X$ be compact complex manifold. A holomorphic vector bundle on $X$ has non-interesting cohomology if $X$ is not compact. The cohomology jump loci of line bundles on a smooth projective curve are clasically called  Brill-Nother loci, see \cite{ACGH}. Cohomology jump loci of vector bundles can be thus considered as generalized Brill-Nother loci.

Let $E$ be a locally free $\cO_X$-coherent sheaf on $X$. By Example \ref{exFIM}, the deformations of $E$ are governed by the {\bf Dolbeault complex} of $E$, 
$$(\Ad(\enmo(E))=(\Gamma(X,\enmo (E)\otimes_{\cO_X}\Omega_X^{0,\bullet}),\bar{\partial}),$$
with the natural DGLA structure. The deformations of $E$ with cohomology constraints  are governed by the DGLA pair
$$(\Ad(\enmo(E)), \Ad(E)).$$

The main issue arising when applying the deformation thechniques of \ref{subDef} - \ref{subDefPairs} is that the DGLA pairs involved are usually infinite-dimensional and not formal, even in the case of smooth projective curves. We will focus here on cases when formality holds.

Consider the moduli space $\cM$ of stable rank $n$ holomorphic vector bundles $E$ with vanishing Chern classes on a compact K\"ahler manifold $X$. We consider the Hodge-theoretic flavored cohomology jump loci 
$$
\sV^{pq}_k(F)\subset \sM
$$
given by
$$
\sV^{pq}_k(F)=\{E\in\sM \mid \dim H^q(X, E\otimes_{\sO_X} F \otimes_{\sO_X} \Omega^p_X)\geq k \}
$$
with the natural scheme structure, for fixed $p$ and fixed poly-stable bundle $F$ with vanishing Chern classes. In this case, the deformation problem with cohomology constraints is governed by the DGLA pair $$(\Ad(\enmo(E)), \Ap(E\otimes F)).$$

\begin{thrm}[\cite{dgms}, \cite{s1}]\label{formal1}
 Let $X$ be a compact K\"ahler manifold.
For any $E\in \sM$, the DGLA pair \\
$$(\Ad(\enmo(E)), \Ap(E\otimes F))$$ is formal. 
\end{thrm}

Let 
\begin{align*}
\sQ(E) &=\{\eta\in H^1(X, \enmo(E))\mid
\eta\wedge\eta=0\in H^2(X, \enmo(E))\}, \\
\sR^{pq}_k(E; F) & =\{\eta\in \sQ(E) \mid  \dim H^q(H^\ubul(X, E\otimes F\otimes \Omega^p_X),\eta \wedge\cdot )\geq k\}.
\end{align*}
The set $\sQ(E)$ is endowed with a natural closed scheme structure defined using the corresponding universal quadratic equations on the affine space $H^1(X, \enmo(E))$. In particular, $\sQ(E)$ has quadratic singularities. The set $\sR^{pq}_k(E; F)$ has a natural closed subscheme structure of $\sQ(E)$ defined using the cohomology jump ideals as in \ref{subEqs} of the  universal complex over $\sQ(E)$ with $\cO_{\sQ(E)}$-linear maps given by the cup-product:
$$
(H^\ubul(X, E\otimes F\otimes \Omega^p_X)\otimes_{\bC}\cO_{\sQ(E)}, \eta_{univ}\wedge\cdot).
$$
Note that this is a complex of free finitely generated $\cO_{\sQ(E)}$-modules such that the entries of the matrices representing the maps in the complex come from linear polynomials in $\bC[H^1(X, \enmo(E))]$. More precisely, $\sQ(E)$ and $\sR^{pq}_k(E;F)$ are the quadratic cone and, respectively, the resonance varieties  as defined in \ref{subEqs} of the cohomology DGLA pair with trivial differentials
\begin{equation}\label{eqCDGA}
(H^\ubul(X,\enmo(E)), H^\ubul(X,E\otimes F\otimes\Omega^p_X)).
\end{equation}

Formality of the DGLA pair implies:

\begin{thrm}\label{thrmHolVB} {\rm (\cite{BW2})} Let $X$ be a compact K\"ahler manifold. Let $E$ and $F$ be a stable and, respectively, a poly-stable holomorphic vector bundle with vanishing Chern classes on $X$. Then there is an isomorphism of formal schemes
$$
 \sV^{pq}_k(F)_{(E)}\cong\sR^{pq}_k(E;F)_{(0)}.
$$ 
\end{thrm}
This isomorphim means that one has explicit simple equations describing locally at $E$ the cohomology jump loci $\sV^{pq}_k(F)$, via the universal equations mentioned above for $\sR^{pq}_k(E; F)$. Note that although the equations are simple, the outcome is not necessarily so.

The isomorphism $$(\sM)_{(E)}\cong \sQ(E)_{(0)},$$ which is a particular case of Theorem \ref{thrmHolVB}, was shown by Nadel \cite{n} and Goldman-Millson \cite{gm}. For rank $n=1$ bundles, the result is due to Green-Lazarsfeld \cite{gl1,gl2}.

By fixing the attention on a generic vector bundle $E$ laying in $\sV^{pq}_k(F)$, that is, if
$$
\dim H^q(X, E\otimes F\otimes \Omega^p_X)=k,
$$
then the generators of the cohomology jump ideals defining $\sR^{pq}_k(E; F)$ are precisely all the entries of the matrices defining the cohomology jump ideals. Since these entries are linear forms, one obtains:

\begin{cor} With the same assumptions as in Theorem \ref{thrmHolVB}, if $$\dim H^q(X, E\otimes F\otimes \Omega^p_X)=k,$$ then $\sV^{pq}_k(F)$ has quadratic algebraic singularities at $E$. 
\end{cor}

This was also showed for $F\otimes\Omega_X^p=\cO_X$ by Martinengo \cite{ma} and B. Wang \cite{w}. 

\begin{rmk}
The assumptions in the Theorem are made such that the controlling DGLA pair is formal and the gauge action of $H^0$ is trivial on $H^1$. As explained in \ref{subEqs}, whenever this is the case, the quadratic cone and the resonance varieties of the cohomology DGLA pair describe locally the moduli space and the cohomology jump loci. We will present repetitions of this type of thinking to obtain similar results in completely different contexts. This is  how one puts the New Principle of \ref{subCC} to work.
\end{rmk}

\section{Deformations of representations and local systems}

To a group $G$, one can associate the set of rank $n$ complex linear representations,
$$
\homo(G,GL(n,\bC)).
$$
When $G$ is finitely presented, this is an affine scheme. Any finite presentation provides equations for this scheme. These are in general hard to compute and to extract information from. We will consider groups $G$ in a more geometric context, namely as fundamental groups of topological spaces.

Let $X$ be a smooth manifold which is of the homotopy type of a finite type CW-complex, and let $x\in X$ be a base point. Then the fundamental group $\pi_1(X,x)$ is finitely presented. We will use the notation 
$$\mathbf{R}(X, n)=\homo(\pi_1(X, x), GL(n, \mathbb{C}))$$
for the moduli space of representations.   

Recall that a (complex) local system of rank $n$ on a topological space is a locally constant sheaf of $\bC$-vector spaces of rank $n$. Every closed point $\rho\in\mathbf{R}(X, n)$ corresponds to a rank $n$ local system $L_\rho$ on $X$. Conversely, a local system $L_\rho$ on $X$ defines a representation $\rho$ uniquely up to conjugation. This corresponds to the fact that local systems do not come with a choice of base point $x$.

Let $W$ be a local system of any rank on $X$. One can  define the cohomology jump loci 
$$
\sV^i_k(W)\subset\bR(X,n)
$$
given by
$$\sV^i_k(W)=\{\rho\in \bR(X, n)\,|\, \dim H^i(X, L_\rho\otimes_\bC W)\geq k\}
$$
with a natural closed scheme structure. For this structure, one uses the cohomology jump ideals of appropriately chosen universal complexes: either formed via a universal twisted cochain complex on the universal cover of $X$, or by using the universal rank $n$ local system on $\bR(X,n)$. Both universal complexes lead to isomorphic cohomology jump loci due to invariance of the latter under quasi-isomorphisms, e.g. \cite{dp, BW2}.

\begin{rmk}\label{rmkPi}
(1) For a group $G$, one can define the cohomology jump loci of $G$ by setting $X=K(G,1)$. 

(2) For a smooth manifold $X$, it is known that the cohomology jump loci $\sV^1_k(\bC_X)$ of rank one representations, which we shall denote by $\Sigma^1_k$ later on, are the same as the those of its fundamental group. In other words, $\Sigma^1_k$ are invariants of $\pi_1(X)$, and can be computed from a finite presentation of $\pi_1(X)$.
\end{rmk}


We address the local structure of the cohomology jump loci $\sV^i_k(W)$ at a representation $\rho\in \bR(X,n)$. Equivalently, we ask how one can describe all the infinitesimal deformations of $\rho$ constrained by the condition  that the degree $i$ cohomology of the corresponding local system $L_\rho$ has dimension $\ge k$.

It is almost true that this infinitesimal deformation problem with cohomology constraints is governed by the DGLA pair given by de Rham complexes
$$(\Ar(\enmo(L_\rho)), \Ar(L_\rho\otimes_\bC W)).$$
Here, for a local system $L$ on $X$,  the de Rham complex $\Ar(L)$ is the complex of global $L$-valued $C^\infty$ forms with the usual differential. Also, $\enmo(L_\rho)$ is the local system of endomorphisms of $L_\rho$. 

To make this completely true, one needs to reintroduce in the picture the choice of base point $x$. Let $$\varepsilon: \Ar(\enmo(L_\rho))\to \mathfrak{g}=\enmo(L_\rho)|_{x}$$ be the restriction map. This is a DGLA map, hence we can also define $$\mathfrak{h}=\varepsilon(H^0(X,\enmo (L_\rho))).$$

Goldman-Millson \cite{gm} attach to an augmented DGLA $(\Ar(\enmo(L_\rho)); \varepsilon)$ a deformation functor $$\Def(\Ar(\enmo(L_\rho)); \varepsilon),$$ defined as the deformation functor of the kernel of $\varepsilon$. They show that the formal scheme of $\bR(X, n)$ at $\rho$ pro-represents this functor. Similarly, one can define cohomology jump loci sub-deformation functors
attached to an augmented DGLA pair $$(\Ar(\enmo(L_\rho)), \Ar(L_\rho\otimes_\bC W); \varepsilon)$$ as in  \ref{subDefPairs}, such that they describe the local structure of $\sV^i_k(W)$ at $\rho$, \cite{BW2}.

\begin{rmk}
In the eye-opening paper \cite{dp}, Dimca and Papadima have addressed the reduced local structure at the trivial representation of cohomology jump loci of representations of fundamental groups. In this case, the de Rham complex $\Ar(X)$ of $X$ plays the controlling role, as it is a module over itself. The results in \cite{BW2} were the outcome of the desire of the authors to extended the results of \cite{dp} beyond the reduced structure and beyond the trivial representation. The solution in terms of DGLA pairs to infinitesimal deformation problems with cohomology constraints seemed natural to us only after reading \cite{dp}.
\end{rmk}

In cases where formality of the (augumented) DGLA pair is present, one can write explicit local equations for $\cV^i_k(W)$ as explained in \ref{subEqs}. Let
$$
\sQ(\rho)=\{\eta\in Z^1(\pi_1(X), \mathfrak{gl}(n, \cc)_{\ad \rho})\,|\,\bar\eta\wedge\bar\eta=0 \in H^2(X, \enmo(L_\rho))\},$$
$$
\sR^i_k(\rho,W)=\{\eta\in  \sQ(\rho)\,|\, \dim H^i(H^\ubul(X, L_\rho\otimes_\bC W), \bar\eta\wedge\cdot)\geq k\},
$$
with the natural scheme structures, where $Z^1$ stands for the vector space of 1-cocycles, $\ad\rho$ is the adjoint representation, and $\bar\eta$ is the image of $\eta$ in cohomology.

Similarly, for a local system $L$ on $X$, define 
$$
\sQ(L)=\{\eta\in H^1(X, \enmo(L)) \mid 
\eta\wedge\eta=0\in H^2(X, \enmo(L))\},
$$
$$\sR^i_k(L,W)=\{ \eta\in \sQ(L)\mid \dim H^i(H^\ubul(X, L\otimes W), \eta\wedge \cdot)\geq k \},
$$
with the natural scheme structures.

\begin{thrm}\label{thmRPP} {\rm (\cite{BW2})}
Let $X$ be a compact K\"ahler manifold,  $\rho\in\mathbf{R}(X,n)$ be a semi-simple representation, and $W$ a semi-simple local system on $X$. Then there is an isomorphism of formal schemes
$$
\sV^i_k(W)_{(\rho)}\cong \sR^i_k(\rho,W)_{(0)} \cong(\sR^i_k(L_\rho,W)\times \mathfrak{g}/\mathfrak{h})_{(0)}.
$$
\end{thrm}

\begin{rmk}
In particular, one has that  $$
\mathbf{R}(X, n)_{(\rho)}\cong \sQ(\rho)_{(0)}\cong (\sQ(L)\times \mathfrak{g}/\mathfrak{h})_{(0)},
$$
a result due to Simpson \cite{s1}.
\end{rmk}

\begin{rmk} With the same assumptions one obtains as before that, if in addition $k=\dim H^i(X,L_\rho)$, then $\sV^i_k(W)$ has quadratic singularities at $\rho$.
\end{rmk}

\begin{rmk} If $\rho$ is a rank one representation, then the affine space $\mathfrak{g}/\mathfrak{h}=0$.  In fact, a more general statement holds as we now describe.
\end{rmk}

Let $\mb(X,n)$ be the moduli space of irreducible local systems on $X$ of rank $n$, see \cite{s3}, and let $W$ another local system. Let $$\Sigma^i_k(W)\subset \mb(X,n)$$ be the cohomology jump loci of irreducible local systems. The controlling DGLA pair in this case is the non-augmented pair $$(\Ar(\enmo(L)), \Ar(L\otimes W)).$$ From the formality of this pair one obtains:

\begin{thrm}\label{thmIrrLS}{\rm (\cite{BW2})} Let $X$ be a compact K\"ahler manifold. Let $L$ be an irreducible local system on $X$, and let $W$ be a semi-simple local system. There is an isomorphism of formal schemes
$$
\mb(X,n)_{(L)}\cong \sQ(L)_{(0)}
$$
inducing an isomorphism
$$
\Sigma^i_k(W)_{(L)}\cong\, \sR^i_k(L,W)_{(0)}. 
$$
\end{thrm}

This has also been proved by a different method by Popa-Schnell \cite{PS} for rank one local systems $L$, $W=\bC_X$, and $X$ a smooth projective complex variety.

\begin{rmk} Due to Simpson's non-abelian Hodge theory \cite{s1}, a similar statement holds for Higgs bundles on smooth projective varieties, via a DGLA pair arising from the Higgs complex, see \cite{BW2}.
\end{rmk}

\section{Line bundles and rank one local systems}

In this section we specialize to line bundles and rank one local systems. For the space of rank one local systems on a topological space $X$ we use the notation
$$
\mb(X) =\Hom (\pi _1(X),\bC^*)=\Hom (H_1(X,\bZ),\bC
^*)=H^1(X,\bC^*).
$$
We will use the following notation for the cohomology jump loci of rank one local systems:
$$\Sigma^i_k=\{L\in\mb(X)\mid \dim H^i(X,L)\ge k\}\subset \mb(X),$$
and identify a representation $\rho\in \Hom (\pi _1(X),\bC^*)$ with the associated local system $L_\rho$. We first recall how the local systems and line bundles are related in the compact case.

\subsection{Local systems vs. Higgs line bundles.} We recall the following from \cite{s1}. Let $X$ be a
compact K\"ahler manifold. Let 
$$
\mdol(X)=\{(E,\theta) \mid {\rm{Higgs\ line\
bundle}}\}=\Pic^\tau(X)\times H^0(X,\Omega_X^1),
$$
where $$\Pic^\tau(X)=\ker \{c_1:\Pic(X)\ra H^2(X,\rr)\}.$$
The space $\mdol(X)$ is endowed with a natural complex structure. There is an isomorphism of real analytic groups
\begin{align*}
\mb(X) &\mathop{\lra}^\sim \mdol(X)\\
 L_\rho &\mapsto (E_\rho,\theta_\rho),
\end{align*}
which can be thought of as the analytification-tropicalization of $\mb(X)$. More precisely, $E_\rho$ is the analytification
$$
E_\rho=L_\rho\otimes_\bC\cO_X,
$$
and $\theta_\rho$ is the tropicalization
$$
\theta_\rho= (1,0)\text{-part of }\log ||\rho|| \in H^1(X,\rr)=\Hom(\pi_1(X),\rr).
$$

Let
$$\cU_B(X)\subset \mb(X)$$ be the unitary local systems. We have $$\cU_B(X)
=\Hom(H_1(X,\bZ), S^1),$$ where $S^1\subset \bC^*$ is the unit
circle, and $$\cU_{Dol}(X)=\Pic^\tau(X).$$ $U_B(X)$ is a totally real
subgroup of $\mb(X)$, whereas $\cU_{Dol}(X)$ is a finite disjoint union of
copies of the Picard variety $\Pic^0(X)$. For example, for a smooth projective curve
of genus $g$, we have $\cU_B(X)=(S^1)^{2g}$ and $\cU_{Dol}(X)={\rm
{Jac}}(X)$, the Jacobian of $X$.

For a Higgs line bundle $(E,\theta)$, let 
$$
H^{pq}(E,\theta) = H^p (H^q(X,E\otimes\Omega_X^\ubul),\theta).
$$
The cohomologies of rank one local systems and of Higgs line bundles are related via:

\begin{thrm} {\rm (\cite{gl1,Be,gl2,s1})} Let $X$ be a compact K\"ahler manifold. Let $L\in\mb(X)$ be a rank one local system with associated Higgs line bundle $(E,\theta)\in\mdol(X)$. Then
$$
H^i(X,L)\cong \bigoplus_{p+q=i}H^{pq}(E,\theta).
$$
\end{thrm}
This has been extended by Simpson beyond the rank one case in \cite{s1,s3}. 

We fix in this section the following notation  about cohomology jump loci:
\begin{align*}
\Sigma^i_k&=\{L\in\mb(X)\mid \dim H^i(X,L)\ge k\}&\subset \mb(X),\\
\sV^{pq}_k &=\{E\in\Pic^\tau(X)\mid \dim H^{q}(X,E\otimes\Omega_X^p)\ge k\}&\subset\Pic^\tau(X),\\
\sV^i_k &=\{E\in\Pic^\tau(X)\mid \dim H^{i}(X,E)\ge k\}&\subset\Pic^\tau(X),\\
\sH^{pq}_k &=\{(E,\theta)\in\mdol(X)\mid \dim H^{p,q}(E,\theta)\ge k \}&\subset \mdol(X).
\end{align*}

\begin{rmk}\label{rmkA} 
Since $\Pic^\tau(X)=\cU_{Dol}(X)\subset\mdol(X)$, one has that
$$
\sV^{pq}_k=\sH^{pq}_k\cap \Pic^\tau(X).
$$
Arapura's trick \cite{a1} shows, via the theorem above, that every irreducible component of a $\sH^{pq}_k$ comes from an irreducible component of $\Sigma^i_k$. 
\end{rmk}

Therefore we can take the point of view that studying  line bundles with torsion first Chern class is a particular case of studying rank one local systems.

\subsection{Local study.} Let $X$ be a complex manifold. Consider the local study of the loci $\Sigma^i_k\subset\mb(X)$ at a rank one representation $\rho$.  The DGLA pair describing the formal schemes $$\Sigma^i_{k,(\rho)}$$ is  $$(\Ar(X), \Ar(L_\rho)),$$
where we denote by $\Ar(X)$ the de Rham complex of $X$. 

Although $\Ar(X)$ might not be formal, from the definition of the deformation functor $\Def (\Ar(X))$ one easily computes that
$$
\mb(X)_{(\rho)}\cong H^1(X,\bC)_{(0)}.
$$
This is of course not necessary, since we have the description
$$
\mb(X)\cong H^1(X,\bC^*).
$$
However, the de Rham complex plays the crucial role in understanding further the properties of the cohomology jump loci.

In particular,  $\Sigma^i_{k,(\bone)}$, the formal schemes at the trivial representation, are described by the DGLA pair consisting of $\Ar(X)$ acting on itself. The Lie structure is trivial here, $\Ar(X)$ is actually a commutative differential graded algebra (CDGA). Using the deformation theory of DGLA pairs as in \ref{subEqs}, one has the following description, due for the reduced structure of $\Sigma^1_k$ to Dimca, Papadima, and Suciu: 

\begin{thrm}\label{thrmDPS}{\rm (\cite{dps})}
Let $X$ be a 1-formal space (e.g. a compact K\"ahler manifold, or the complement of a hyperplane arrangement), that is, the de Rham complex is 1-formal. Then $$\Sigma^1_{k,(\bone)}\cong \cR^1_{k,(0)},$$
and $\sR^1_k$ is the tangent cone at $\bone\in\mb(X)$ of $\Sigma^1_k$. In particular, the irreducible components of $\Sigma^1_k$ passing though $\bone$ are subtori.
\end{thrm}
Here,  
$
\cR^i_k
$  
are the cohomology jump loci of the universal Aomoto complex
$$
(H^\ubul (X,\bC)\otimes_\bC\bC[H^1(X,\bC)], w_{univ}\cup .).
$$
Set-theoretically,
$$
\sR^i_{k}=\{w\in H^1(X,\bC)\mid \dim H^i((H^\ubul (X,\bC), w\cup .)\ge k\}.
$$
This follows from making explicit  for this case the section \ref{subEqs}.

\begin{rmk}\label{rmkHAI} (Hyperplane arrangements I.)
This result, for the reduced germs of  $\Sigma^i_{k}$ at the trivial representation, for the case when $X$ is the complement of a hyperplane arrangement, had been obtained earlier by Esnault, Schechtman, and Viehweg \cite{ESV}. By Orlik-Solomon \cite{OS}, the cohomology ring of $X$ is a combinatorial invariant of the hyperplane arrangement, that is, it only depends on the lattice of intersection of the hyperplanes and not on their positions. Hence, the components through $\bone$ of  $\Sigma^i_k$ can be detected combinatorially in this case from the resonance varieties $\cR^i_k$. In the context of hyperplane arrangements, the resonance varieties have been introduced by Falk, see the survey articles by Suciu referred to in the Introduction.
\end{rmk}

\begin{rmk} (Hyperplane arrangements II.) A very difficult folklore conjecture states that all the cohomology jump loci $\Sigma^i_k$ of the complement of a hyperplane arrangement are combinatorial invariants. A particular case of this conjecture states that the eigenvalues of the monodromy on the Milnor fiber of a hyperplane arrangement are determined combinatorially.  The best result so far in this direction is the recent paper of Papadima and Suciu \cite{PaSu}. Surprisingly, they use the cohomology jump loci of rank two local systems and resonance varieties defined over finite fields to address the rank one case. 
\end{rmk}

\begin{rmk}\label{rmkHt}
The de Rham complex of $X$ encodes the homotopy type of the topological space $X$.  A consequence is that the formal scheme $$\Sigma^i_k(X)_{(\bone)}$$ at the trivial representation only depends on the $i$-homotopy type of $X$. For the underlying reduced germs, this is due to Dimca-Papadima \cite{dp}. The cohomology jump deformation functors from \ref{subDefPairs} give the statement for the non-reduced structure as well.
\end{rmk}

\begin{rmk}\label{rmkGys} (Gysin model)
Typically $\Ar(X)$ is not formal. Hence the explicit description of $\Sigma^i_k(X)_{(\bone)}$ depends on finding at least a finite model for $\Ar(X)$, that is, a finite dimensional CDGA quasi-isomorphic with $\Ar(X)$. The crucial observation of Dimca-Papadima \cite{dp} is that, if $X$ is a smooth quasi-projective variety, there is always such a model. This is the {\it Gysin model}  of Morgan \cite{Mo}, defined with respect to any smooth projective compactification $$j:X\ra Y$$ with a simple normal crossings divisor $Y\setminus X=\bigcup_{i\in S}D_i$ as boundary. Here $D_i$ are the irreducible components. For $I\subset S$, let $D_{I}=\cap_{i\in I}D_i$. The Gysin CDGA $\cA^\ubul$ is defined by setting 
$$
\cA^i =\bigoplus_{p+q=i} \cA^{p,q},\quad\quad \cA^{p,q}=\bigoplus_{|I|=q}H^p(D_I,\bC).
$$
The multiplication is the cup-product
$$
\cA^{p,q}\otimes \cA^{p',q'}\mathop{\lra}^\cup \cA^{p+p',q+q'}.
$$
The differential on $\cA^\ubul$ is defined using the Gysin maps from intersections of the boundary divisors and it satisfies
$$
\cA^{p,q}\lra \cA^{p+2,q-1}.
$$
The Gysin model $\cA^\ubul$ is in fact the (total complex of the) $E_2$ page of the Leray spectral sequence for the open embedding $j$:
$$
\cA^{p,q}=E_2^{p,q}=H^p(Y,R^qj_*\bC_X) \Rightarrow \bH^i(Y, Rj_*\bC_X)=H^i(X,\bC).
$$
It is well-known that this spectral sequence degenerates at the $E_3$ page, \cite{PeS}. In other words, the cohomology CDGA of $\cA^\ubul$ is the cohomology of $X$. The more general statement, showing by \ref{subDefPairs} that the infinitesimal deformations of the trivial rank one local system on $X$ with cohomology constraints are governed by the Gysin model, is the following:
\end{rmk}

\begin{thrm} {\rm \cite{Mo}} Let $X$ be a smooth complex quasi-projective variety. Then there is a quasi-isomorphism of CDGAs between the de Rham complex and the Gysin model,
$$
\Ar(X) \cong \cA^\ubul.
$$
\end{thrm}

Hence, as explained in \ref{subEqs}, the resonance varieties of the Gysin model describe the local structure of $\Sigma^i_k$ at the trivial representation.

Using the positive grading on the Gysin model, Dimca and Papadima draw global consequences from local information, generalizing Theorem \ref{thrmDPS}:

\begin{thrm}{\rm (\cite{dp})}\label{thrmDP}
Let $X$ be a smooth quasi-projective complex variety. Then the components of $\Sigma^i_k$ containing the trivial representation are subtori of $\mb(X)$.
\end{thrm}

Further improvements of this result are described next.

\subsection{Structure of the cohomology jump loci.} We describe now in more detail the global structural results on cohomology jump loci of rank one local systems and line bundles.  The proofs of the structure theorems for cohomology jump loci are all similar in plan:  a combination of local study, information from mixed Hodge structure theory, and arithmetic.

\begin{thrm}\label{thrmStr1} {\rm (\cite{a2})}
Let $X$ be the complement in a compact K\"ahler manifold $\bar{X}$ of a divisor with simple normal crossings singularities. Assume that $H^1(\bar{X},\bC)=0$. Then $\Sigma^i_k$ are finite union of unitary translates of subtori of $\mb(X)$. The same holds for $\Sigma^1_1$ with no assumption on $\bar{X}$, with the extra feature that the positive dimensional components are torsion translates.
\end{thrm}

\begin{thrm}\label{thrmStr2}\label{thrmBW1} {\rm (\cite{BW1})} Let $X$ be a smooth complex quasi-projective variety. Then the cohomology jump loci of rank one local systems $\Sigma^i_k$ are a finite unions of torsion translates of  subtori of $\mb(X)$. 
\end{thrm}

\begin{theorem}\label{thrmStr3}{\rm (\cite{W13})} Let $X$ is a compact K\"ahler manifold.  
Then the cohomology jump loci of line bundles $\Sigma^i_k$ and $\sV^{pq}_k$ are a finite unions of torsion translates of  subtori of $\mb(X)$ and, respectively, $\Pic^\tau(X)$. 
\end{theorem}

Here, by subtori of $\mb(X)$ we mean {\it affine tori} $(\bC^*)^s$. By subtori of $\Pic^\tau(X)$ we mean {\it compact complex subtori}, that is, connected compact complex Lie subgroups; their underlying topological spaces are real tori. 

By Remark \ref{rmkA}, the statement about $\sV^{pq}_k$ in the third theorem follows from the statement about $\Sigma^i_k$.

\begin{rmk} The theorems build upon a long list of partial results. The chronology is the following:

- Beauville \cite{Be} showed the linearity property  of $\sV^1_1$ for a compact K\"ahler manifold and conjectured the third theorem for $\sV^{pq}_k$;

- Green-Lazarsfeld \cite{gl2} showed the linearity of $\sV^{pq}_k$ and $\Sigma^1_k$ for a compact K\"ahler manifold;

- Arapura \cite{a1} showed that $\Sigma^i_k$ are finite unions of unitary translated subtori for a compact K\"ahler manifold. He also showed that $\sH^{pq}_k$ are finite unions of unitary translated subtori of $\Pic^\tau(X)$ times vector subspaces of $H^0(X,\Omega_X^1)$, recovering thus the result of \cite{gl2}.

- Simpson \cite{s2}  proved the theorems for $X$ smooth projective, namely by showing that the components of $\Sigma^i_k$, and hence of $\sV^{pq}_k$, are torsion translated subtori; 

- Arapura \cite{a2} showed the first theorem. A small gap in the proof was filled later in \cite{Art-all}.

- Campana \cite{Ca} showed the third theorem for $\Sigma^1_k$;

- Pink and Roessler \cite{PR} reproved the statement about $\sV^{pq}_k$ for a smooth projective variety via reduction modulo $p$;

- in \cite{Bu-ULS} we showed that $\Sigma^i_k\cap \cU_{B}(X)$, the cohomology jump loci of unitary rank one local systems on a smooth quasi-projective variety $X$, and their Hodge theoretic refinements, are finite unions of torsion translated subtori of $\Pic^\tau(X)$ times rational convex polytopes. Combined with Theorem \ref{thrmStr1}, this gives Theorem \ref{thrmStr2} under the assumptions on the compactification as in Theorem \ref{thrmStr1}, as remarked by Dimca and Libgober.

- Dimca-Papadima \cite{dp} showed Theorem \ref{thrmDP}.

- The final statements of the last two theorems were after this just highly non-trivial statements about existence of torsion points in each irreducible component of $\Sigma^i_k$, this being shown recently in \cite{BW1} and, respectively, \cite{W13}.

\end{rmk}

The three theorems suggest that the following question might have a positive answer:

\begin{que}{\it 
Let $X$ be the complement of a simple normal crossings divisor in a compact K\"ahler manifold. Are the cohomology jump loci $\Sigma^i_k$ always finite unions of torsion translated subtori?}
\end{que}

\subsection{Fibrations}

The structure theorem is more refined in some cases. For a compact K\"ahler manifold $X$, let $$a_X:X\ra Alb(X)$$ be the Albanese map. For a compact analytic variety $Y$, we say that that $Y$ has {\it maximal Albanese dimension} if $\dim Y$ equals the dimension of the Albanese of any smooth model of $Y$.

\begin{thrm}\label{thrmGL2} {\rm(\cite{gl2})} Let $X$ be a compact K\"ahler manifold. Let $Z$ be a positive dimensional irreducible component of $\sV^{i}_k\subset\Pic^\tau(X)$. Then there exists an analytic dominant map
$$f:X\ra Y$$
with connected fibers to an analytic variety $Y$ of maximal Albanese dimension and
$$
\dim Y\le i,
$$
such that
$Z$ is a translate of $f^*\Pic^0(Y)$.
\end{thrm}

As an application, a more careful analysis of the reasons behind Theorem \ref{thrmGL2}  allowed Pareschi and Popa to generalize the classical Castelnuovo - de Franchis theorem from surfaces to other dimensions:

\begin{thrm}{\rm (\cite{PP})}
Let $X$ be a compact K\"ahler manifold of maximal Albanese dimension and irregularity $q(X)$. If 
$$
\chi(\omega_X)<q(X)-\dim X
$$
then $X$ admits a surjective morphism with connected fibers onto a normal compact analytic variety $Y$ with
$$ 0 < \dim Y < \dim X$$ and with any smooth model $\tilde{Y}$ of Albanese general type and with $\chi(\omega_{\tilde{Y}})>0$.
\end{thrm}

In the non-compact case, one has the following:

\begin{thrm}\label{thrmA2}{\rm (\cite{a2})} 
Let $X$ be the complement in a compact K\"ahler manifold $\bar{X}$ of a divisor with simple normal crossings singularities. Assume that $H^1(\bar{X},\bC)=0$. Let $W$ be a positive dimension irreducible component of  $\Sigma^i_k$. Then there exists an analytic map
$$
f:X\ra T
$$
to an extension $T$ of a compact complex torus by an affine torus, such that $W$ is a  translate of $f^*\mb(T)$. 
\end{thrm}

In  Theorems \ref{thrmGL2} and \ref{thrmA2}, the translates are by torsion elements whenever the assumptions of Theorem \ref{thrmStr1} - \ref{thrmStr3} hold.

The components of the support loci $\Sigma^1_1$ of the first cohomology are completely understood in terms of {\it fibrations}, i.e. surjective analytic maps with connected general fiber, onto curves. Generalizing previous results of \cite{gl1, Be, Cat, gl2, a1} from the compact case to the non-compact case, Arapura proved:

\begin{thrm}{\rm (\cite{a2})} Let $X$ be the complement in a compact K\"ahler manifold of a divisor with simple normal crossings singularities. Then there is a bijection between:
\begin{itemize}
\item the set equivalence classes of fibrations
$$
f:X\ra C,
$$
where $C$ is a smooth Riemann surfaces minus finitely many points with $\chi(C)<0$, under the equivalence relation $f\sim f'$ if there exists an isomorphism $g:C\ra C'$ with $g\circ f=f'$, and
\item the irreducible positive-dimensional components of $\Sigma^1_1$ containing the trivial local system.
\end{itemize}
To a fibration $f$ corresponds the component $f^*\mb(C)=f^*\Sigma^1_1(C)$.
\end{thrm} 

\begin{rmk}
This result has recently been generalized to normal complex varieties by Arapura, Dimca, and Hain \cite{ADH}. In this paper, they show that $\Sigma^1_k$ for a normal variety behaves as for a smooth variety. This opens up the question of how much of what is known about cohomology jump loci of smooth varieties goes through to varieties with singularities.
\end{rmk}

\begin{rmk} {\rm (\cite{a2})}
An interesting application is that the cardinality $N_b(X)$ of the set of equivalence classes of fibrations onto curves as above with fixed Betti number $b_1(C)=b$ depends only on the fundamental group $\pi_1(X)$. Indeed, it is known that the support of $\Sigma^1_k$ can be computed from $\pi_1(X)$. By the theorem, $N_b(X)$ is the number of irreducible components of $\Sigma^1_1$ of dimension $b$.
\end{rmk}

By the torsion translate property, the components of $\Sigma^1_1$ not passing through the trivial local system, must contain a torsion point. Hence, one can account for these components by applying the theorem to a finite cover of $X$. Or, equivalently, one can use the orbifold language, to which we refer to \cite{Art-all}. Moreover, in \cite{Dcons, Dirr, Art-all}, all the positive dimensional components components of $\Sigma^1_k$ for $k>1$ are also accounted for. The structure of $\Sigma^1_k$ forces constraints on fundamental groups $\pi_1(X)$, in light of Remark \ref{rmkPi}. Here are some of the additional properties:

\begin{thrm} Let $X$ be a smooth quasi-projective complex variety. Then:
\begin{itemize}
\item {\rm (\cite{Dirr})} Any irreducible component of $\Sigma^1_k$ is also one of $\Sigma^1_{l}$ for any $1\le l\le k$.
\item {\rm (\cite{DPSa})}Any two distinct irreducible components of $\Sigma^1_k$ intersect in at most finitely many points, and the intersection points are torsion.
\item {\rm (\cite{Art-all})} If $Z$ and $Z'$ are distinct irreducible components of $\Sigma^1_k$ and, respectively, $\Sigma^1_{k'}$, and if $\rho\in Z\cap Z'$, then $\rho\in \Sigma^1_{k+k'}$.
\end{itemize}
\end{thrm}

\subsection{Dimension bounds}

Theorem \ref{thrmGL2} implies a dimension bound:

\begin{thrm}\label{thrmGLbd} {\rm (\cite{gl1,gl2})} Let $X$ be a compact K\"ahler manifold.  The cohomology support loci $\sV^i_1\subset\Pic^\tau(X)$ satisfy
$$
\codim \sV^i_1\ge \dim a_X(X) - i.
$$
\end{thrm}

\begin{rmk}
This codimension bound can also be thought of as a special case of the codimension bound (\ref{eqN}) for DGLA pairs. Indeed, consider $\sV^i_1$ locally at trivial line bundle $\cO_X$. This is described by the DGLA pair (\ref{eqCDGA}) which simplifies in this case to
$$
(H^\ubul(X,\cO_X),H^\ubul(X,\cO_X)).
$$
Let $\ed(X)$ be the exactness degree of this pair as defined in \ref{subDB}. Here the quadratic cone of the pair is the affine space $$\sQ=H^1(X,\cO_X).$$ The universal Aomoto complex in this case is
$$
(H^\ubul(X,\cO_X)\otimes_\bC\cO_{\sQ},w_{univ}).
$$
Applying (\ref{eqN}), one has that 
$$
\codim \sV^i_1\ge \ed(X) - i.
$$
Then the theorem follows after identifying 
$$
\ed(X)=\dim a_X(X),
$$
this equality being explicitly stated in \cite[Theorem A]{LP}.
\end{rmk}

The bound on the codimension of $\sV^i_1$ implies the generic vanishing:

\begin{thrm}\label{thrmGVGL}{\rm (\cite{gl1})} Let $X$ be a compact K\"ahler manifold. Then for a general line bunde $L\in\Pic^\tau(X)$, 
$$
H^i(X, L)=0\quad\quad\text{for}\quad\quad i<\dim a_X(X).
$$
\end{thrm}

Theorem \ref{thrmGLbd} has been generalized to Hodge pieces and to local systems  by Popa and Schnell:

\begin{thrm}\label{thrmPS} {\rm (\cite{PS})} Let $X$ be a smooth complex projective variety of dimension $n$. 
\begin{itemize}
\item The cohomology support loci $\sV^{pq}_1\subset\Pic^\tau(X)$ satisfy
$$
\codim \sV^{pq}_1\ge |p+q-n|-\delta(a_X), 
$$
\item The cohomology support loci $\Sigma^i_1\subset\mb(X)$ satisfy
$$
\codim \Sigma^i_1\ge 2(|i-n|-\delta(a_X)).
$$
\end{itemize}
\end{thrm}

The first part in the theorem implies the second part by the correspondence between local systems and Higgs line bundles. Here $\delta(a_X)$ is the defect of semi-smallness of the Albanese map, 
$$
\delta(a_X)=\max_{l\in\nn}(2l-\dim X+\dim A_l),
$$
where 
$$
A_l=\{y\in Alb(X)\mid \dim f^{-1}(y)\ge l\}.
$$
A weaker statement for abelian varieties was also proved by Kr\"amer-Weissauer \cite{KW}.

\begin{rmk}
We like to think of Theorem \ref{thrmPS} as, again, a special case of the codimension bound (\ref{eqN}). In this case, the controlling DGLA pair for $\sV^{pq}$ locally at the trivial line bundle is
$$
(H^\ubul(X,\cO_X),H^\ubul(X,\Omega_X^p)),
$$
according to (\ref{eqCDGA}). The quadratic cone of this pair is $\sQ=H^1(X,\cO_X)$ again. The universal Aomoto complex of this pair is
$$
(H^\ubul(X,\Omega_X^p)\otimes\cO_{\sQ}, w_{univ}).
$$
Let $\ed(X,p)$ be the exactness degree of this complex. By (\ref{eqN}),
$$
\codim \sV^{pq}_1\ge \ed(X,p)-q.
$$
By Hodge symmetry, it is enough to consider the case $p+q\le n$. Then the result follows from the bound
$$
\ed(X,p)\ge n-p-\delta(a_X),
$$
shown implicitly in \cite{PS}. 
\end{rmk}

The proof in \cite{PS} of Theorem \ref{thrmPS}  goes in two steps. The first step is to identify a rather natural class of coherent sheaves $\cF$ on an abelian variety $A$ of dimension $n$ such that the relative cohomology jump loci 
$$
\sV^i_k(\cF)=\{L\in\Pic^0(A)\mid \dim H^i(A,\cF\otimes L)\ge k\}\subset \Pic^0(A)
$$
satisfy automatically the bound
$$
\codim \sV^i_1(\cF)\ge i.
$$
These are called {\it GV sheaves}. This bound is by duality related to the bound (\ref{eqN}) obtained by classical homological algebra. The second step is to note that the derived direct image of the trivial line bundle under the Albanese map of $X$ decomposes into special sheaves as above by the Decomposition Theorem for proper maps. The special sheaves are:

\begin{thrm}{\rm (\cite{PS})} Let $A$ be a complex abelian variety. Let $M$ be a mixed Hodge module on $A$ with underlying filtered $\cD$-module $(\cM,F)$. Then
$$
\codim \sV^i_1(\gr^F_k\cM)\ge i
$$
for all $k$ and $i$.
\end{thrm}

\subsection{Application to abelian varieties and compact complex tori.} {\rm (\cite{PPS})} For a stronger statement and a generalization of the previous theorem to compact complex tori, see the recent article of Pareschi, Popa, and Schnell \cite{PPS}. The GV sheaves, and their related cousins the M-regular sheaves, that is, coherent sheaves $\cF$ satisfying the stronger inequality
$$
\codim \sV^i_1(\cF)>i,
$$
have been used  to prove the following numerical characterizations in terms of plurigenera:

\begin{thrm} {\quad}

\begin{itemize}

\item {\rm (\cite{CH})} A smooth complex projective variety $X$ is birational to an abelian variety if and only if $P_1(X) = P_2(X) = 1$ and $q(X) = \dim X$ .

\item {\rm (\cite{PPS})} A compact K\"ahler manifold $X$ is bimeromorphic to a compact complex torus if and only if $\dim H^1(X,\bC) = 2\dim X$ and $P_1(X) = P_2(X) = 1.$
\end{itemize}
\end{thrm}

\subsection{Constructible complexes on tori.} Let $\cK$ be a $\bC$-constructible complex on a complex manifold. Consider the relative hypercohomology jump loci of rank one local systems
$$
\Sigma^i_k(\cK)\subset \mb(X)
$$
given by
$$
\Sigma^i_k(\cK)=\{L\in \mb(X)\mid \dim \bH^i(X,\cK^\ubul\otimes _\bC L)\ge k\}.
$$

\begin{thrm}\label{thrmGaLo} {\rm (\cite{GaLo})}  Let $\cK$ be a perverse sheaf on the affine torus $X=(\bC^*)^n$.
Then $
\Sigma^i_k(\cK)$ is a finite unions of translates of subtori of $\mb(X)$ and
$$
\codim \Sigma^i_1(\cK)\ge i.
$$
\end{thrm}

We expect the translates to be by torsion elements.

As a corollary, one has the following generic vanishing, the perverse sheaf counterpart of Theorem \ref{thrmGVGL}.

\begin{thrm} 
Let $X\subset (\bC^*)^n$ be a closed subvariety of dimension $d$ such that the shifted constant sheaf $\bC_X[d]$ is a perverse sheaf on $(\bC^*)^n$ (e.g. a locally complete intersection). Then
$$
(-1)^d\chi(X)\ge 0.
$$
\end{thrm}

\begin{rmk} This result had appeared already in \cite{LS} and was attributed to Laumon in \cite{GaLo}. The positivity of the signed Euler characteristic has been later conjectured for any closed subvariety of $(\bC^*)^n$ in \cite{HS}. Counterexamples to this conjecture were recently constructed in \cite{BW14}.
\end{rmk}

In the compact complex torus case, one has:

\begin{thrm}\label{thrmPPS} {\rm (\cite{W13, PPS})} Let $X$ be a compact complex torus. Let $\cK$ be a perverse sheaf on $X$ coming from a polarizable real Hodge module.   
Then $
\Sigma^i_k(\cK)$ is a finite unions of torsion translates of subtori of $\mb(X)$.
\end{thrm}

Here a polarizable real Hodge module is a generalization of unitary local systems, as opposed to the usual polarizable Hodge modules being generalizations of quasi-unipotent local systems of geometric origin. The result of \cite{W13} deals with polarizable Hodge modules, the extension to real polarizable Hodge modules is due to \cite{PPS}.

\section{Fundamental groups and homotopy type}
\subsection{Fundamental groups}
Let $X$ be a connected smooth manifold. As pointed out in Remark \ref{rmkPi}, the moduli space of rank one local systems $\mb(X)$ and the cohomology jump loci $\Sigma^1_k$ are determined by the fundamental group $\pi_1(X)$. The statement that each irreducible component of $\Sigma^1_k$ is a torsion translated subtorus puts nontrivial conditions on the fundamental group of $X$. 

First, let us recall some notations in group theory. Given a group $G$, we denote its derived subgroup $[G,G]$ by $DG$ or $D^1G$. Let $D^nG=D(D^{n-1}G)$. $\rank(G/D^nG)$ is defined to be $$\sum_{1\leq j\leq n}\dim_\qq(D^{j-1}G/D^jG)\otimes_\zz \qq.$$ $\rank(G/D^nG)$ can be a natural number or $+\infty$.  

A group $G$ is called {\it virtually nilpotent} if it has a finite index subgroup $G'$ which is nilpotent. A group $G$ is called {\it polycyclic} if there is a sequence of subgroups $G=G_0\rhd G_1\rhd \cdots \rhd G_n=\{0\}$ such that each $G_j$ is a normal subgroup of $G_{j-1}$ and each $G_{j-1}/G_{j}$ is cyclic. 

The following theorem is proved using Simpson's result about $\Sigma^i_k$ of smooth projective varieties. 

\begin{thrm}{\rm (\cite{AN})}
Let $X$ be a smooth projective variety, and let $G=\pi_1(X)$. Suppose that $DG$ is finitely generated and $G/D^nG$ is solvable of finite rank for some $n$. Then there are normal subgroups $P\supseteq Q\supseteq D^nG$ so that
\begin{enumerate}
\item $G/P$ is finite,
\item $P/Q$ is nilpotent,
\item $Q/D^nG$ is a torsion group.
\end{enumerate}
\end{thrm}

When $G$ is polycyclic, $D^nG=0$ for sufficiently large $n$. Moreover the assumption of the theorem is automatically satisfied. Thus we have the following consequence. 
\begin{cor}{\rm (\cite{AN})}
Let $X$ be a smooth projective variety, and let $G=\pi_1(X)$. Suppose $G$ is polycyclic. Then $G$ is virtually nilpotent. 
\end{cor}

By extending Simpson's result to $\Sigma^1_k$ of compact K\"ahler manifolds, in \cite{Ca} Campana generalized the above results to the case when $X$ is a compact K\"ahler manifold. 

Let now $G_n=[G_{n-1},G]$ denote the {\it lower central series} of a group $G$. The {\it Chen groups} of a group $G$ are the successive quotients in the lower central series of $[G,D^2G].$ A remarkable result of Cohen-Schenck \cite{CS}, conjectured by Suciu, is that the resonance varieties $\cR^1_1$, as in Remark \ref{rmkHAI}, determine combinatorially the ranks of the Chen groups of the fundamental group of a complement of hyperplane arrangements:

\begin{thrm} {\rm (\cite{CS})} Let $\theta_k(G)$ be the rank of the $k$-th Chen group of the fundamental group $G$ of a complement $U$ of a hyperplane arrangement. Let $h_k$ be the number of $k$-dimensional irreducible components  of the resonance variety
$$
\cR^1_1(U)=\{ w\in H^1(U,\bC)\mid H^1(H^\ubul (U,\bC),w\cup .)\ne 0\}.
$$
Then
$$
\theta_k(G)=(k-1)\sum_{m\ge 2}h_m\binom{m+k-2}{k}.
$$
\end{thrm}

This is a consequence of 1-formality, reducedness of the resonance varieties, and some other special properties that fundamental groups of complements of hyperplane arrangements enjoy. Beyond such groups, it remains to be determined what is the precise relation between the possibly-nonreduced cohomology jump loci and Chen groups.

\subsection{Homotopy type of smooth projective varieties and compact K\"ahler manifolds}
As we have discussed in Theorem \ref{thrmStr1}, Theorem \ref{thrmStr2} and Theorem \ref{thrmStr3} and thereafter, the special geometry of smooth quasi-projective varieties and compact K\"ahler manifolds has strong implications on the cohomology jump loci of rank one local systems $\Sigma^i_k$. Next, we will show some examples of how to use cohomology jump loci to detect the non-algebraicity of compact K\"ahler manifolds. 

Every one dimensional compact complex manifold is projective hence K\"ahler. In dimension two, there are compact complex manifolds that are not K\"ahler, for example the Hopf surface or more generally any compact complex surface $X$ with odd $b_1(X)$. There are also compact K\"ahler surfaces that are not projective, for example a general 2-dimensional complex torus. However, using the classification of complex surfaces, Kodaira proved that every compact K\"ahler surface can be deformed to a smooth projective surface. In higher dimensions, Kodaira problem asks whether a higher dimensional compact K\"ahler manifold can always be deformed to a smooth projective variety. 

In \cite{v}, Voisin gave a negative answer to Kodaira problem by a topological method. More precisely, Voisin constructed compact K\"ahler manifolds that are not of the same homotopy type as any smooth projective variety. She also showed how to modify the example to obtain compact K\"ahler manifolds that are not of the same rational (or real) homotopy type as any smooth projective variety. Let us review Voisin's construction first. 

Given a compact complex torus $T$ of dimension at least two. Suppose there is a torus endomorphism $\phi$ of $T$ whose induced endomorphism $H^1(\phi)$ on $H^1(T, \zz)$ satisfies some special properties. Voisin proved that the special endomorphism forces that the rational Neron-Severi group $NS_{\qq}(T)=H^{1,1}(T, \cc)\cap H^2(T, \qq)$ is trivial. Hence $T$ is not an abelian variety. 

Fix a compact complex torus $T$ with the special endomorphism $\phi$ as above. Voisin constructed a 4-dimensional compact K\"ahler manifold $X$ by blowing up a sequence of points and surfaces on $T\times T$. Voisin managed to encode the information of the endomorphism $\phi$ into the cohomology ring $H^\ubul(X, \zz)$. She showed that if another compact K\"ahler manifold $Y$ has the same cohomology ring, i.e., $H^\ubul(Y, \zz)\cong H^\ubul(X, \zz)$ as graded commutative rings, then $Alb(Y)$ is not an abelian variety, and hence $Y$ is not projective. 

Voisin also gave other constructions so she can encode the information of $\phi$ in the $\rr$-coefficient cohomology ring $H^\ubul(X, \rr)$. Thus the optimal result of \cite{v} is the following. 

\begin{thrm}{\rm \cite{v}}
There exists a compact K\"ahler manifold $X$ which is not of the real homotopy type of any smooth projective variety. 
\end{thrm}

In \cite{W13}, Voisin's theorem is reproved using cohomology jump loci. A compact K\"ahler manifold $X$ is constructed which is very similar to Voisin's example. Now, the information of the endomorphism $\phi$ is now encoded in the cohomology jump locus $\Sigma^2_1(X)$ of rank one local systems. Thus, by Voisin's argument using rational Neron-Severi group, we obtain the following theorem. 

\begin{thrm}{\rm \cite{W13}}
There exists a compact K\"ahler manifold $X$ with the following property. Suppose $Y$ is another compact K\"ahler manifold with $H^1(Y, \cc)\cong H^1(X, \cc)$, or equivalently there is an isomorphism of germs $\mb(Y)_{(0)}\cong \mb(X)_{(0)}$. If the isomorphism $\mb(Y)_{(\bone)}\cong \mb(X)_{(\bone)}$ induces an isomorphism $\Sigma^2_1(Y)_{(\bone)}\cong \Sigma^2_1(X)_{(\bone)}$, then $Y$ is not projective. 
\end{thrm}

As mentioned in Remark \ref{rmkHt}, the pairs of germs $\Sigma^i_k(X)_{(\bone)}\subset \mb(X)_{(\bone)}$ are determined by the real $i$-homotopy type of $X$ for all $i$ and $k$. Thus we can slightly improve Voisin's result to the following. 

\begin{cor}
There exists a compact K\"ahler manifold $X$ that is not of real 2-homotopy type of any smooth projective variety. 
\end{cor}

\section{Small ball complements and Bernstein-Sato ideals}
We now let $X=\bC^n$ and let $$f:X\ra \bC$$ be a holomorphic function. For a point $x\in f^{-1}(0)$, let 
$$U_{f,x}=B_x\cap (X\setminus f^{-1}(0)),$$
where $B_x$ is a small ball centered at $x$ in $X$.

\begin{conj} Let $U_{f,x}$ be the complement of a hypersurface singularity germ as above. Then 
the cohomology jump loci of rank one local systems $\Sigma^i_k\subset\mb(U_{f,x})$ are finite unions of torsion translates of subtori.
\end{conj}

\begin{rmk}
It is known that the small ball complements $U_{f,x}$ are 1-formal spaces, see \cite{DH}. Therefore Theorem \ref{thrmDPS} applies, and hence the conjecture is true for those irreducible components of $\Sigma^1_k(U_{f,x})$ containing the trivial rank one local system, that is, those components are subtori.
\end{rmk}

\begin{rmk} This conjecture is stated erroneously as a theorem in \cite{Li-gap}. We have also believed it to be proved in {\it loc. cit.} until the writing of \cite{BW2} when we realized that the rather difficult Theorem \ref{thrmStr2} should in principle be easier to prove. A red flag about \cite{Li-gap} was raised then. Scouting for the details, we remarked the following. Let $L$ be a rank one local system on a a smooth quasi-projective variety $X$. It is known that if $L$ is unitary, then the cohomology $H^\ubul(X,L)$ admits a mixed Hodge structure \cite{t}. The attached theory of weights leads, via a well-known argument of Deligne, to a technique for degeneration of associated spectral sequences, cf.  Remark \ref{rmkGys}. In a more generalized context, this is also the reason why Theorem \ref{thrmPPS} works for polarizable real Hodge modules. However, if $L$ is not unitary, it is not clear how to construct mixed Hodge structures nor how to attach weights. This leaves unproved beyond the unitary case Lemma 3.5 of \cite{Li-gap} on the degeneration of a certain spectral sequence. At least for this reason, the main result of \cite{Li-gap} is unfortunately still open, and we stated this as the above conjecture.
\end{rmk}

The above conjecture generalizes a classical result in singularity theory. Recall that, by Milnor, the small ball complement $U_{f,x}$ is a smooth locally trivial fibration over a small punctured disc $\Delta^*$ in $\bC$ via the map $f$. The homotopy class of the fiber is called {\it the Minor fiber of $f$ at $x$} and it is denote by $F_{f,x}$. The conjecture implies the well-known:

\begin{thrm}{\bf (Monodromy Theorem)}
Let $f: (\bC^n,x)\ra (\bC,0)$ be the germ of a holomorphic function. Then the eigenvalues of the monodromy on the cohomology $H^i(F_{f,x},\bC)$ of the Milnor fiber are roots of unity for every $i$.
\end{thrm}

The precise relation between the conjecture and the Monodromy Theorem is the following. Assume, for simplicity, that there is locally at $x$ a splitting
$$
f=\prod_{1\le i\le r} f_i
$$
such that 
\begin{equation}\label{eqF}
\text{the }f_i\text{ define mutually distinct, reduced, irreducible analytic branches of }f.
\end{equation}
In this case,
$$
\mb(U_{f,x})=(\bC^*)^r,
$$
and a rank one local system $L$ on $U_{f,x}$ is identified with the $r$-tuple $\lam\in(\bC^*)^r$ given by the monodromies around each of the branches.

\begin{lemma} Let $f$ and $x$ be as above, satisfying (\ref{eqF}). Let 
$$
diag : \bC^*\lra \mb(U_{f,x})=(\bC^*)^r
$$
be the diagonal embedding $\gamma\mapsto (\gamma,\ldots,\gamma)$. Then intersection of the cohomology support loci of rank one local systems
$$
\bigcup_i\Sigma^i_1(U_{f,x})\subset \mb(U_{f,x})
$$
with the image of the diagonal embedding is the set of eigenvalues of the Milnor monodromy on the cohomology of $F_{f,x}$.
\end{lemma}

\begin{rmk} When (\ref{eqF}) does not hold, one has a similar statement after adjusting the definition of the diagonal embedding. The lemma follows immediately from \cite[Prop. 3.31, Thm. 4]{Bu-BS}. Note that, as pointed out by Liu-Maxim \cite{LM}, in all the statements in \cite{Bu-BS} where the uniform support $Supp^{unif}\psi_F\bC_X$ of the Sabbah specialization complex appears, the {\it unif} should be dropped to conform to what is proven in \cite{Bu-BS}. There are two other different proofs of this lemma known to us: one using the Leray spectral sequence for the Milnor fibration, the other one using perverse sheaves related to the open embedding of $U_{f,x}$ in $B_x$. 
\end{rmk}

A different approach to the Monodromy Theorem was given by Malgrange and Kashiwara. For a holomorphic function $f:X\ra \bC$, one defines the {\it Bernstein-Sato polynomial} $b_f(s)$ to be the monic generator of the ideal of polynomials $b(s)\in\bC[s]$ satisfying
$$
b(s)f^s=Pf^{s+1},
$$
for some operator $$P\in\cO_X[\pa_{1}, \ldots, \pa_n, s],$$ where $\cO_X$ are the holomorphic functions on $X$, and $\pa_i={\pa}/{\pa x_i}$ are the usual derivations with respect to a choice of coordinates on $X=\bC^n$. A non-trivial result of Bernstein and Sato guarantees the existence of a non-zero $b_f(s)$ for any $f$. If $f$ is a polynomial, then $\cO_X$ can be replaced with the polynomial ring in $n$ variables.

One has the classical main result about Bernstein-Sato polynomials, letting $\Exp(\al)=\exp (2\pi i \al)$:

\begin{thrm} Let $f:\bC^n\ra \bC$ be a holomorphic function. Then:
\begin{itemize}
\item {\rm (\cite{Ka})}The roots of $b_f(s)$ are negative rational numbers.
\item {\rm (\cite{Mal, Ka2})}The set
$$
\{\Exp (\alpha)\mid \al \text{ is a root of }b_f(s)\}\subset \bC^*
$$
is the set of all eigenvalues of the Milnor monodromy on the cohomologies of Milnor fibers $F_{f,x}$ of $f$ at points $x\in f^{-1}(0)$.
\end{itemize}
\end{thrm}

We have taken in \cite{Bu-BS} the point of view that this result of Malgrange and Kashiwara hides a more general statement about rank one local systems and we tried to uncover what this statement is. Our answer is in terms of the {\it Bernstein-Sato ideals} $$B_{F,x} \subset \bC[s_1,\ldots, s_r]$$
attached to the collection of functions
$$
F=(f_1,\ldots, f_r),
$$
with, for simplicity, $f=\prod_{i=1}^rf_i$ as in (\ref{eqF}). The ideal $B_{F,x}$ is defined as the the ideal generated by the polynomials $b\in\bC[s_1,\ldots, s_r]$ such that
$$
b(s_1,\ldots,s_r)\prod_{i=1}^rf_i^{s_i}=P\prod_{i=1}^rf_i^{s_i+1}
$$
for some operator 
$$
P\in \cO_{X,x}[\pa_1,\ldots,\pa_n,s_1,\ldots,s_r],
$$
where $x_i$ are local coordinates at $x$.

For $y\in f^{-1}(0)$ close to $x$, one has the cohomology support loci $\Sigma^i_1(U_{f,y})$ in $\mb(U_{f,y})$. Since now there might be branches $f_{i}$ with $f_{i}(y)\ne 0$, it is possible that $\mb(U_{f,y})$ is a smaller affine torus than $(\bC^*)^r$. Fixing the coordinates $t_1,\ldots ,t_r$ on $(\bC^*)^r$, the smaller affine torus is cut out by the equations $t_{i}=1$ for those $i$ with $f_i(y)\ne 0$. By eliminating these equations, we define the {\it uniform cohomology support loci}
$$
\Sigma^i_1(U_{f,y})^{unif}\subset (\bC^*)^r,
$$
all sitting inside the same affine torus for all $y$ close to $x$. Consider the zero locus $$Z(B_{F,x})\subset\bC^r$$ of the Bernstein-Sato ideal.

\begin{thrm} {\rm (\cite{Bu-BS})} With the notation as above,
$$
\Exp (Z(B_{F,x}))\supset\mathop{\bigcup_{i\in\bZ}}_{y\in f^{-1}(0)\text{ close to }x} \Sigma^i_1(U_{f,y})^{unif},
$$
and we conjectured that equality holds.
\end{thrm}

\end{document}